\documentclass[12pt]{amsart}
\usepackage{amssymb,latexsym,amsmath,amscd}
\theoremstyle{plain}
\newtheorem{theorem}{Theorem}[section]
\newtheorem{proposition}[theorem]{Proposition}
\newtheorem{corollary}[theorem]{Corollary}

\theoremstyle{definition}
\newtheorem{definition}[theorem]{Definition}
\newtheorem{remark}[theorem]{Remark}
\newtheorem{notation}[theorem]{Notation}
\newtheorem{example}[theorem]{Example}

\newcommand{\ra}{\rightarrow}
\newcommand{\lra}{\longrightarrow}

\newcommand{\PP}{\mathbf{P}}

\newcommand{\NN}{\mathbf{N}}
\newcommand{\ZZ}{\mathbf{Z}}

\newcommand{\bin}[2]{ {{#1} \choose {#2}} }
\newcommand{\OO}{\mathcal{O}}

\newcommand{\FF}{\mathcal{F}}

\newcommand{\tr}{\rm{tr}}

\newcommand{\degr}{\rm{deg }}

\newcommand{\re}{{\rm{reg}}}

\title[Hilbert Functions of Irreducible Arithmetically Gorenstein Schemes]
{Hilbert Functions of Irreducible Arithmetically Gorenstein
Schemes}

\author{Nero Budur}
\address{Department of Mathematics \\Johns Hopkins University \hfil\break\indent 3400 North Charles Street\\
Baltimore, MD 21218, USA} \email{nbudur@math.jhu.edu}

\author{Marta Casanellas}
\address{Departament d'Algebra i Geometria, Facultat de
Matem\`atiques \hfil\break\indent Gran Via 585, 08071 Barcelona, Spain}
\email{casanell@mat.ub.es}

\author{Elisa Gorla}
\address{Department of Mathematics
\\ University of Notre Dame \hfil\break\indent 255 Hurley Hall, Notre Dame, IN 46556-4618, USA}
\email{egorla@nd.edu}

\thanks{Research of the third author partially supported by ``Istituto Italiano di Alta Matematica
Francesco Severi.'' The three authors were also supported by the
organization of Pragmatic 2000. The research in the paper was done
in June 2000 at the University of Catania - Italy.}

\begin{document}

\maketitle

\begin{abstract}
 In this paper we compute the Hilbert functions of
irreducible (or smooth) and reduced arithmetically Gorenstein
schemes that are twisted anti-canonical divisors on arithmetically
Cohen-Macaulay schemes. We also prove some folklore results
characterizing the Hilbert functions of irreducible standard
determinantal schemes, and we use them to produce a new class of
functions that occur as Hilbert functions of irreducible (or
smooth) and reduced arithmetically Gorenstein schemes in any
codimension.
\end{abstract}

\section*{Introduction}

There is a simple characterization of the functions that arise as
Hilbert functions of arithmetically Cohen-Macaulay schemes.
Nevertheless, very little is known about the Hilbert functions of
{\em irreducible} arithmetically Cohen-Macaulay schemes. Harris
proved that the $h$-vector of an irreducible aCM scheme of
positive dimension is the same as the $h$-vector of a zero-scheme
satisfying the Uniform Position Property (see \cite{Ha}). However,
no characterization is available yet for these $h$-vectors. In
codimension 2, the question has been completely answered following
a different approach. Notice that in codimension 2, due to the
Hilbert-Burch theorem, arithmetically Cohen-Macaulay schemes and
standard determinantal schemes coincide. A standard determinantal
scheme is defined by the (homogeneous) maximal minors of a
suitable homogeneous matrix of forms, see Definition \ref{stdet}.
In the first part of the paper, we state some folklore facts about
irreducible standard determinantal schemes, leaving the
proofs for the Appendix. First, we express the Hilbert function of
any standard determinantal scheme in terms of its degree matrix
(Proposition \ref{hilbstd}). Then, we characterize the Hilbert
functions of irreducible and reduced standard determinantal
schemes in $\PP^n$ of any codimension in terms of the entries of
the degree matrix (Theorem \ref{cor-st-det}). In particular, since
any standard determinantal scheme is an arithmetically
Cohen-Macaulay scheme, we obtain a large class of numerical
functions that occur as the Hilbert functions of some irreducible
and reduced arithmetically Cohen-Macaulay schemes.

If one restricts attention to {\it irreducible} arithmetically
Gorenstein schemes, the question of characterizing their Hilbert
functions has been answered only in the codimension 3 case, thanks
to the Buchsbaum-Eisenbud structure theorem (see \cite{DNV}).
However, the question is still open in higher codimension. In the
second part of the paper, we produce a new class of Hilbert
functions that occur for irreducible and reduced (respectively, irreducible
and smooth) arithmetically Gorenstein schemes of any
codimension (Corollary \ref{last.thm}, and Corollary
\ref{smooth}). The strategy is to show (Theorem
\ref{lastprop}) that, for an aCM subscheme $S\subset\PP^n$ which
is Gorenstein in codimension one, a general element of the linear
system $|mH-K|$ determines an irreducible arithmetically
Gorenstein divisor whose $h$-vector can be written in terms of the
Hilbert function of $S$. Here $H$ is a hyperplane section of $S$ by a
hyperplane that meets it properly, $K$ a
canonical divisor, and $m\gg 0$
(Theorem 3.2 also contains an estimate of how big $m$ can be chosen).
The Corollaries mentioned above
are obtained by combining this result with the folklore results of
the first part.

In the first section we recall a few facts about Hilbert
functions. In the second section we state the folklore facts
mentioned above. In the third section we draw our main
conclusions. The Appendix contains the proofs omitted in Section 2.

\textbf{Acknowledgment}: We thank everybody in Pragmatic 2000,
especially A. Bigatti, A.V. Geramita, J. Migliore, C. Peterson and
A. Ragusa. We also thank the referee for comments which greatly
improved on the first version.

\section{Preliminaries}

Let $S$ be a closed subscheme of the projective $n$-space $\PP
^n=\PP^n(k)$, where $k$ is an algebraically closed field. Let
$I_S$ be the saturated homogeneous ideal corresponding to $S$ in
the polynomial ring $R=k[x_0,\ldots,x_n]$.

The numerical
function  \[\begin{array}{rcl} H_S:\NN & \longrightarrow & \NN \\
n & \longmapsto & \dim_k (R/I_S)_n \end{array}\] is called the
{\it Hilbert function} of $S$. The formal series  \[P_S
(z)=\sum_{n\geq 0} H_S (n)z^n\] is called the {\it Hilbert series}
of $S$. It is well-known that the Hilbert series of $S$ can be
expressed in the rational form
\[ P_S(z)=\frac{h_S(z)}{(1-z)^{d+1}}\ ,\]
where $h_S(z)$ is a polynomial with integer coefficients such that
$h_S(1)=\deg S$ is the degree of the scheme, and $d$ is the
dimension of $S$. The polynomial $$h_S(z)=\sum_{i=0}^s h_i z^i,$$
with $h_s\ne 0$, is called the {\it h-polynomial} of $S$, and the
vector $(h_0,\ldots,h_s)$ defined by the coefficients of $h_S(z)$
is called the {\it h-vector} of $S$. We will say that an
$h$-vector $h=(h_0,\ldots,h_s)$ has length $s$.

By an {\it arithmetically Cohen-Macaulay} (abbreviated aCM)
projective scheme we mean a projective scheme whose coordinate
ring is Cohen-Macaulay. For an aCM scheme $S$,
set $H_S(-1)=0$, and define
$$\begin{array}{rcl}
\Delta^1 H_S (t) & = & H_S(t)-H_S(t-1),\\
\Delta^r H_S (t) & = & \Delta^{r-1} H_S (t)-\Delta^{r-1} H_S (t-1).
\end{array}$$
$\Delta^r H_S$ is called the {\it $r$-th difference} of $H_S(t)$ and it is the
Hilbert function of the $r$-th general hyperplane section $X$ of $S$.

The problem of characterizing those numerical functions that occur
as Hilbert functions of  schemes with given properties has been
studied extensively. There is a simple characterization of the
Hilbert series of aCM projective schemes (see e.g. \cite{V},
Theorem 1.5). However, in general very little is known about the
Hilbert series of irreducible aCM schemes. The case of irreducible
aCM schemes of codimension 2 is better understood, thanks to the
structure theorem of Hilbert and Burch (see \cite{BH}, Theorem
1.4.17). The numerical functions that can occur as Hilbert
functions for reduced, irreducible aCM schemes of codimension 2
are characterized in \cite{GM} and in \cite{V}, Theorem 2.3. There
is no analogous characterization in higher codimension.

By an {\it arithmetically Gorenstein} (abbreviated aG) scheme we
mean a projective scheme whose coordinate ring is Gorenstein. In
particular, any aG scheme is an aCM scheme. Even in the case of
arithmetically Gorenstein projective schemes there is no complete
characterization of the Hilbert series. A necessary, but not
sufficient, condition for a polynomial to be the $h$-polynomial of
some aG scheme is symmetry in the coefficients, i.e. if
$h(z)=1+h_1 z+\ldots+h_{s-1}z^{s-1}+h_s z^s$, then $h_s=1$ and
$h_i=h_{s-i}$ for all $i=1,\ldots,s-1$.

For codimension three aG
schemes, the Hilbert series can be characterized using  the
structure theorem of Buchsbaum-Eisenbud (see \cite{V}, Theorem 2.18).
For irreducible codimension three aG schemes, a
characterization of the Hilbert series is given by De Negri and
Valla in \cite{DNV}. Our Corollary \ref{last.thm} gives new
examples in each codimension greater than three of Hilbert series
occurring for irreducible aG schemes. Corollary \ref{smooth} does
the same, under the extra assumption that the irreducible aG schemes
be smooth.

Finally, we briefly recall the definition of Castelnuovo-Mumford
regularity.
\begin{definition}\label{defreg} A coherent sheaf $\FF$ on $\PP^n$ is said to be
$m$-\emph{regular} if
$$H^i(\PP^n, \FF(m-i))=0$$
for all $i>0$. The \emph{regularity} or \emph{Castelnuovo-Mumford
regularity} of $\FF$ is
$${\re} {(\FF)}={\text{min}} \{m | \FF {\rm{\ is\ }}m \rm{-regular}\}.$$
\end{definition}

By a theorem of Serre (\cite{H}, Theorem 3.5.2), any $\FF$ is
$m$-regular for some $m$. It follows from the Castelnuovo-Mumford
Theorem that if $\FF$ is an $m$-regular coherent sheaf on $\PP^n$,
then $\FF(k)$ is generated as an $\OO_{\PP^n}$-module by its
global sections for all $k\geq m$ (see e.g. \cite{M}, Theorem 1.1.5).

If $V \subset \PP^n$ is a subscheme we define the
\emph{Castelnuovo-Mumford regularity} of $V$ as the
Castelnuovo-Mumford regularity of its ideal sheaf,
$\re {(V)}:=\re {(\mathcal{I}_V)}$.

\begin{remark}\label{reg}
(i) It is well known (and easy to prove) that if
$$0 \lra \dots \lra \oplus_{j} R(-a_{i,j}) \lra \dots \lra
\oplus_{j} R(-a_{0,j}) \ra I_V \lra 0$$ is a graded minimal free
$R$-resolution of $I_V$, then
$${\re} {(\mathcal{I}_V)}=\max_{i,j} \{a_{i,j}-i \}$$
(see for instance \cite{M}, Remark 1.1.6).\\
 (ii) The
Castelnuovo-Mumford regularity of a scheme $S \subset \PP^n$
should not be confused with the {\it regularity index} of a
scheme, $r(S)$, that is the minimum degree in which the Hilbert
function of $S$ agrees with the Hilbert polynomial (see \cite{BH}
for more details). If $S \subset \PP^n$ is an aCM scheme of
dimension $d$ , then $r(S)=\re {(S)}-d-1$ (this follows from
\cite{BH}, Theorem 4.4.3 (b)).
\end{remark}
%ooooooooooooooooooooooooooooooo

\section{Hilbert functions of standard determinantal
schemes}\label{folklore}

The results stated in this section are regarded as folklore.
We will compute the Hilbert function of a standard determinantal scheme in terms
of the degree matrix associated to it. In particular, we will derive formulas for the
degree and the Castelnuovo-Mumford regularity of the scheme.
Complete proofs using standard methods are given in the Appendix.

\begin{definition}\label{stdet} A subscheme $S \subset \PP^n$ is called {\it standard
determinantal} if $I_S$ is generated by the maximal minors of an
$l\times (l+c-1)$ homogeneous matrix $M=(g_{ij})$ representing a
morphism $$\phi:F=\bigoplus_{j=1}^{l+c-1} R(a_j) \lra
G=\bigoplus_{j=1}^l R(b_j)$$ of free graded $R$-modules.
Here $c$ is the codimension of $S$ and we assume that $a_1 \leq
\dots \leq a_{l+c-1}$ and $b_1 \leq \dots \leq b_l$. The degree
matrix of $M$ is the matrix
 $U=(u_{ij})$  whose entries are the degrees of the
entries of $M$.  We will call {\em degree matrix} any matrix of
integers that is the degree matrix associated to some homogeneous
matrix of polynomials.
\end{definition}

\begin{remark}\label{deg} (i) In the notation above, the entries of $U$ increase
from right to left and
from top to bottom: $u_{i,j} \geq u_{k,r}$ if $i \geq k$ and $j
\leq r$. \\
(ii) We will assume without loss of generality that the degree
matrix has the property $g_{ij}=0$ if $u_{ij}\le 0$, and $\degr$
$g_{ij}=u_{ij}$ if $u_{ij}>0$. Note that $g_{ij}$ could be $0$
even if $u_{ij}>0$. \\
(iii) Since  $U=(u_{ij})$ is the degree matrix of a homogeneous
matrix, one has
\begin{equation}\label{pp}\sum_{v=1}^{s}
u_{i_v,j_v}=\sum _{v=1}^{s} u_{i_v, j_{\pi (v)}},\end{equation}
for every permutation $\pi$ of $\{1,\ldots ,s\}$.
\end{remark}

Any standard determinantal scheme is an aCM scheme and its minimal
free resolution is given by the Eagon-Northcott complex (see proof
of Proposition \ref{hilbstd} in the Appendix). Moreover, in
codimension 2, any aCM scheme is a standard determinantal scheme
due to the Hilbert-Burch Theorem.

Next, we show how to recover the Hilbert function of a standard determinantal
scheme from its degree matrix.

\begin{notation} Let $U=(u_{ij})$ be a matrix of size $l\times (l+c-1)$.
By an $m\times m$ {\it submatrix} $V$ of $U$ we mean the $m\times
m$ matrix obtained from the elements of $U$ given a choice of $m$
rows and $m$ columns. For a square matrix $V$ we say that $\dim
V=m$ if $V$ is a $m\times m$ matrix. By $(V_1|V_2|\ldots )\subset
U$ we mean a choice of square submatrices $V_j$ of $U$ with $\dim
V_1=l$, such that: the choice of columns for $V_{j+1}$ is strictly
to the right-hand side of any column chosen for $V_j$, and the
choice of rows for $V_{j+1}$ is made from the choice of rows for
$V_j$ (so $0\leq \dim V_{j+1} \leq \dim V_j \leq l$). Denote by
$\tr (V_1|V_2|\ldots)$ the sum of the traces of $V_j$ (see Example \ref{excomp}).
\end{notation}

With the notation above, we have the following result whose
proof is left to the Appendix. We
thank the referee who pointed out a similar formula obtained by Ausina
and Ballesteros in the unpublished paper \cite{AB}
(see also Section 5 of \cite{AB2}).
Our results were obtained independently,
and we give different formulas than those of Ausina
and Ballesteros. The main tool used here is the Eagon-Northcott resolution,
as in \cite{AB}.

\begin{proposition}\label{hilbstd}
Let $S \subset \PP^n$ be a standard determinantal scheme of
codimension $c$ with degree matrix $U=(u_{i,j})$, $i=1, \dots,l$,
$j=1, \dots, l+c-1$. Then
$$H_S(t) ={{t+n}\choose{n}}+ \sum_{(V_1|V_2|\ldots)\subset
U}(-1)^{1+\dim V_2+\ldots}\bin{t+n-\tr(V_1|V_2|\ldots)}{n},$$
$$\degr(S)=1+\frac{1}{c!}\sum_{(V_1|V_2|\ldots)\subset U}
(-1)^{c+1+\dim
V_2+\ldots}(\tr(V_1|V_2|\ldots)-1)\cdot\ldots\cdot({\tr(V_1|V_2|\ldots)}-c).$$
The Castelnuovo-Mumford regularity of $S$ is
$$\re (S)={\tr(V_1^0 | V_2^0| \dots|V_c^0)}-c+1,$$
where $V_1^0$ is the submatrix of $U$ formed by the first $l$
columns, and $V_{k+1}^0=(u_{l,l+k})$ for $k=1, \dots c-1$ (i.e.
$V_{k+1}^0$ is the $(l+k)$-th element in the last row of $U$).
\end{proposition}

\begin{remark}
(i) The expression for $\deg S$ in Proposition \ref{hilbstd} is a
new version of Porteus' formula in the case of standard
determinantal schemes (see \cite{ACGH}, II.4.2). The advantage is
that this formula involves only the entries of the degree matrix
of $S$, while Porteus' formula involves Chern classes.\\
 (ii) From Proposition \ref{hilbstd} and  Remark \ref{reg} we have that the index of
regularity of $S$ is
$$r(S)={\tr(V_1^0 | V_2^0| \dots|V_c^0)}-n.$$
\end{remark}

\begin{example}\label{excomp}
Let $S \subset \PP^5$ be a standard determinantal scheme of
codimension $3$ with degree matrix
$$U=\left( \begin{array}{cccc} 2 & 2 & 2 & 1 \\
3& 3 & 3 & 2 \end{array} \right)$$ To compute the degree and the
regularity of $S$ we have to consider the following combinations
of submatrices of $U:$
$$
\left(\begin{array} {cc} 2 & 2 \\ 3& 3
\end{array}  \begin{array}{cc}
\cdot & \cdot \\ \cdot & \cdot \end{array} \right), \,
\left(\begin{array} {cc} 2 & 2 \\ 3& 3
\end{array} \vline \begin{array}{cc}
2 & \cdot \\ \cdot & \cdot \end{array} \right), \,
\left(\begin{array}{cc} 2 & 2 \\ 3& 3
\end{array} \vline \begin{array}{cc}
2 \; \vline & 1 \\ \cdot & \cdot \end{array} \right), \,
\left(\begin{array} {cc} 2 & 2 \\ 3& 3
\end{array} \vline \begin{array}{cc}
\cdot & \cdot \\ 3 & \cdot \end{array} \right), \,
\left(\begin{array} {cc} 2 & 2 \\ 3& 3
\end{array} \vline \begin{array}{cc}
\cdot & \cdot \\ 3 \; \vline & 2 \end{array} \right),
$$
$$
\left(\begin{array}{cc} 2 & 2 \\ 3& 3
\end{array} \vline \begin{array}{cc}
2 & 1 \\ 3 & 2 \end{array} \right), \, \left(\begin{array} {cc} 2&
2\\ 3& 3
\end{array} \vline \begin{array}{cc}
\cdot & 1 \\ \cdot & \cdot \end{array} \right), \,
\left(\begin{array} {cc} 2 & 2 \\ 3& 3
\end{array} \vline \begin{array}{cc}
\cdot & \cdot \\ \cdot & 2 \end{array} \right), \,
\left(\begin{array} {cc} 2 & \cdot \\ 3& \cdot
\end{array} \begin{array}{cc}
2 & \cdot \\ 3 & \cdot \end{array} \right), \,
\left(\begin{array}{ccc} 2 & \cdot & 2 \\ 3& \cdot & 3
\end{array} \vline \begin{array}{r}
 1 \\
\cdot \end{array} \right), \,
$$
$$
\left(\begin{array}{ccc} 2 & \cdot &
2 \\ 3& \cdot & 3
\end{array} \vline \begin{array}{r}
 \cdot \\
2 \end{array} \right), \, \left(\begin{array}{cccc} 2 & \cdot &
\cdot & 1 \\ 3 & \cdot & \cdot & 2
\end{array}\right), \,
\left(\begin{array}{cccc} \cdot & 2 & 2 & \cdot  \\ \cdot & 3 & 3
& \cdot \end{array}\right), \,
\left(\begin{array}{ccc}  \cdot & 2
& 2
\\ \cdot & 3 & 3
\end{array} \vline \begin{array}{r}
 1 \\
\cdot \end{array} \right), \,
\left(\begin{array}{ccc}  \cdot & 2&
2
\\ \cdot & 3 & 3
\end{array} \vline \begin{array}{r}
\cdot \\
2 \end{array} \right), \,
$$
$$ \left(\begin{array}{cccc} \cdot& 2 & \cdot & 1 \\ \cdot & 3 &
\cdot & 2
\end{array}\right), \,
\left(\begin{array}{cccc}  \cdot & \cdot & 2 & 1 \\  \cdot & \cdot
& 3 & 2
\end{array}\right) .
$$
The degree of $S$ is
\begin{multline*} \degr(S)=1+\frac{1}{6}\left( 4\cdot 3\cdot
2-6\cdot 5\cdot 4+7\cdot 6\cdot 5-7\cdot 6\cdot 5 +9\cdot 8\cdot
7+8\cdot 7\cdot 6-5\cdot 4\cdot 3-6\cdot 5\cdot 4 \right. \\
\left. +4\cdot 3\cdot 2-5\cdot 4\cdot 3-6\cdot 5\cdot 4+ 3\cdot
2\cdot 1+4\cdot 3\cdot 2-5\cdot 4\cdot 3-6\cdot 5\cdot 4+3\cdot
2\cdot 1+3\cdot 2\cdot 1 \right)=  46.
\end{multline*}
For the degree matrix $U$ we have $$(V_1^0|V_2^0|V_3^0)=\left( \begin{array}{cc}2 & 2 \\
3 & 3
\end{array} \vline \begin{array}{cc}
\cdot & \cdot \\ 3 \; \vline & 2
\end{array} \right),$$
hence Proposition \ref{hilbstd} tells us that the
Castelnuovo-Mumford regularity of $S$ is $\re (S)=10-3+1=8.$
\end{example}
\vskip 5mm

We focus now on determining the numerical functions that are
Hilbert functions of standard determinantal schemes with
interesting properties such as reduced and irreducible.

Next, given a matrix satisfying certain positivity conditions on
the entries, we construct a reduced standard determinantal scheme
that has this degree matrix. The existence of such a scheme also
follows from work of Trung (\cite{T}). A self-contained proof of the
following proposition is given in the Appendix.

\begin{proposition}\label{reduced}
Let $U=(u_{i,j})$ be a degree matrix of size $l\times (l+c-1)$,
satisfying the condition (ii) of Remark \ref{deg}. Suppose that
$u_{i,i+c-1}>0$, $i=1, \dots, l$. Then, for any $n$ with $n\geq
c\geq 1$, there exists a reduced standard determinantal scheme $X
\subset \PP^n$ of codimension $c$ with degree matrix $U$.
\end{proposition}

We can now obtain a large class of Hilbert functions that occur as
Hilbert functions of reduced and irreducible, arithmetically Cohen-Macaulay
schemes.

\begin{theorem}\label{cor-st-det} Let $U=(u_{i,j})$ be a degree matrix
of size $l\times (l+c-1)$. Let $n$ be an integer, $n > c\geq 1$. Then
$$H(t) ={{t+n}\choose{n}}+ \sum_{(V_1|V_2|\ldots)\subset
U}(-1)^{1+\dim V_2+\ldots}\bin{t+n-\tr(V_1|V_2|\ldots)}{n}$$ is
the Hilbert function of a non-degenerate, irreducible and reduced standard
determinantal scheme $X \subset \PP^n$ if and only
if $u_{i,i+c}>0$ for $i=1,\dots ,l-1$.
\end{theorem}

\begin{proof}
It follows from a theorem of Trung (see \cite{T}, Theorem 4.4)
that if $S\subset \PP ^n$ is a standard determinantal scheme of
codimension $c\ge 2$ with degree matrix $U$, then $S$ is the
hyperplane section of a normal, standard determinantal
reduced, irreducible scheme $S'\subset \PP ^{n+1}$ of codimension
$c$ by a hyperplane that meets it properly if and only if $u_{i,i+c}>0$ for $i=1,\dots ,l-1$. Clearly,
$S$ and $S'$ will have the same degree matrix. Then the theorem
follows directly from Proposition~\ref{hilbstd}.
\end{proof}

%ooooooooooooooooooooooooooooooo

\section{Hilbert functions of irreducible arithmetically Gorenstein
schemes}

In this section we obtain a set of functions that are Hilbert
functions of a large class of irreducible and reduced
arithmetically Gorenstein schemes.
In order to be able to use Bertini's Theorem (as in Theorem~\ref{lastprop}),
we will be working over an algebraically closed field $k$ of characteristic
$0$.

Recall that a noetherian ring $A$ (respectively, a noetherian scheme $X$)
satisfies the condition $G_r$,
{\it Gorenstein in codimension less or equal $r$}, if
every localization $A_P$ at a prime ideal $P\subset A$ (respectively, every
local ring $\OO_{X,x}$)
of codimension less then or equal to $r$ is a Gorenstein local ring. In other words,
the non locally-Gorenstein locus
has codimension greater than $r$ (see \cite{Har2} for more details).

Let $S\subset \PP^n$ be a codimension $c$ scheme satisfying
property $G_1$, and let $I_S$ be its saturated homogeneous ideal.
A {\it divisor} $D$ on $S$ is a
generalized divisor in the sense of \cite{Har2}.

We will denote a dualizing sheaf on $S$ by
$\omega_S$ and the corresponding canonical divisor by $K$. We will
denote the canonical module of $S$ by $K_S$, that is
$$K_S=Ext^c_R(R/I_S,R)(-n-1) \ $$
which is isomorphic to $H^0_{\ast}(\omega_S)$. Here,
$H^i_{\ast}(\mathcal{F}):=\oplus_{t \in \ZZ}H^i(\PP ^n,
\mathcal{F}(t))$ for any sheaf $\mathcal{F}$. Finally, we will denote by $\omega_S^\vee$
the $\OO_S$-dual of the canonical sheaf $\omega_S$.

In order to construct arithmetically Gorenstein schemes we will
use the following result:
\begin{proposition}\label{KMMNP}(Corollary 5.5, \cite{KMMNP})
Let $S\subset \PP ^n$ be an aCM subscheme satisfying $G_1$, $K$ a canonical divisor
on $S$, and $H$ the hyperplane section. Then every element of the linear system
$|mH-K|$ is arithmetically Gorenstein.
\end{proposition}

Next, we compute the Hilbert functions of these arithmetically
Gorenstein schemes, in order to obtain numerical functions that occur
as  Hilbert functions of irreducible arithmetically Gorenstein
schemes. We express the Hilbert functions of these aG divisors on a scheme $S$
only in terms of
the Hilbert function of $S$ and its regularity.
Notice that all the $h$-vectors that arise this way are of
\emph{decreasing type}, in the sense that if
$\Delta h_Y(j_0)<0$ for some $j_0$ then $\Delta h_Y(j)<0$ for all $j\geq j_0$.

While it has been proved that the $h$-vector of an arithmetically Gorenstein,
reduced and irreducible
scheme of codimension 3 is of decreasing type (see \cite{DNV}), it is not known
whether the same holds in codimension 4 or higher.

\begin{theorem}\label{lastprop}
Let $S \subset \PP ^n$ be an irreducible and reduced aCM scheme of
dimension $d\ge 2$, satisfying property $G_1$. Denote by
$\Delta^dH_S$  the Hilbert function of the $d^{\rm{th}}$-general
hyperplane section of $S$. Set
$$r:={\min}\ \{ i | \Delta^d H_S(i)=\deg S \}=\re (S)-1.$$
Let $Y$ be a general element in the linear system $|mH-K|$ for
$m\ge\max\{2r-d,\re (\omega_S^{\vee})\}$. Then $Y$ is an irreducible and
reduced arithmetically Gorenstein scheme whose $h$-vector is of
decreasing type and satisfies
$$h_Y(t)= \left\{ \begin{array}{cr}
\Delta^dH_S(t), & t \leq r \\
\deg S, & r \leq t \leq m-r+d \\
\Delta^dH_S(m-t+d), & t\geq m-r+d.
\end{array} \right.$$
\end{theorem}

\begin{proof}
Let $Y \subset S$ be a general element of the linear system
$|mH-K|$, $m \in \ZZ$. For $m\ge \re (\omega_S^{\vee})$, this linear
system is base point free and, by Bertini's Theorem (see \cite{J},
p. 89), the general element $Y$ is irreducible. Moreover, $Y$ is
arithmetically Gorenstein by Proposition \ref{KMMNP}. Let $I_S$
and $I_Y$ be the saturated homogeneous ideals of $S$ and $Y$ as
subschemes of $\PP ^n$ and let $\mathcal{I}_{Y,S}$ be the
sheafification of the ideal $I_Y/I_S \subset R/I_S$. We have
$$\mathcal{I}_{Y,S} \cong \mathcal{O}_S(K-mH) \cong \omega_S(-m) \
$$ and $S$ is aCM of dimension $>1$. Therefore
$H^1_{\ast}(\mathcal{I}_S)=0$, and we obtain the exact sequence
$$ 0 \ra I_S \ra I_Y \ra H^0_{\ast}(\omega_S)(-m) \ra 0 \ $$
by taking cohomology in $0 \ra \mathcal{I}_S \ra \mathcal{I}_Y \ra
\mathcal{I}_{Y,S} \ra 0$. Thus we get the following equality
on Hilbert functions for every $t$:
\begin{equation}\label{Hf} \ H_Y(t) = H_S(t)-H_{K_S}(-m+t) \ .\end{equation}

Since $S$ is aCM, the dual of the
resolution of $S$ is a resolution for $K_S(n+1)$, see
\cite{M}, Remark 1.4.8. Hence one can write the Hilbert function of the canonical module in terms of the
$h$-vector of $S$:
$$\Delta^{d+1}H_{K_S}(t)=h_S(d+1-t).$$
Therefore by (\ref{Hf}) we get
$$\Delta^{d+1}H_Y(t)=\Delta^{d+1}H_S(t)-\Delta^{d+1}H_{K_S}(-m+t)=h_S(t)-h_S(d+1-t+m),$$
and for any integer $t>0$:
$$h_Y(t)=\Delta^dH_Y(t)=\sum_{i=0}^t h_S(i) - \sum_{i=0}^t h_S(d+1-i+m)=$$
$$=\Delta^dH_S(t)+ \sum_{j=0}^{m-t+d}h_S(j)- \sum_{j=0}^{m-t+d}h_S(j)-
\sum_{j=m-t+d+1}^{m+d+1}h_S(j)=$$
$$=\Delta^dH_S(t)+\Delta^dH_S(m-t+d)-\Delta^dH_S(m+d+1) \ .$$
Here $\Delta^dH_S$ is the Hilbert function of the
$d^{\rm{th}}$-general hyperplane section of $S$, which is a set of
points.  Set $r:=\min \{ i | \Delta^d H_S(i)={\degr} S \}$. Since
$m \geq r-d$, we obtain that
$$h_Y(t)= \Delta^dH_S(t)+ \Delta^dH_S(m-t+d)-{\degr} S.$$
If $m \geq 2r-d$, then the $h$-vector of $Y$ is:
$$h_Y(t)= \left\{ \begin{array}{cr}
\Delta^dH_S(t), & t \leq r \\
{\degr} S, & r \leq t \leq m-r+d \\
\Delta^dH_S(m-t+d), & t\geq m-r+d
\end{array} \right.$$
Notice that $r=\re (S)-1$. This follows from Remark \ref{reg} and
from the fact that $r=r(S)+d$, where $r(S)$ is the index of
regularity of $S$.

To prove that $h_Y$ is of decreasing type, we will see that if
$\Delta h_Y(j_0)<0$ for some $j_0$, then $\Delta h_Y(j) <0$ for
any $j \geq j_0$. For $t \leq r$ we have $\Delta
h_Y(t)=\Delta^{d+1}H_S(t)$ which is strictly positive because it is
the $h$-vector of an aCM scheme. For $r
\leq t \leq m-r+d$ we have $\Delta h_Y(t)=0$, and for any $t \geq
m-r+d$ we have $\Delta
h_Y(t)=\Delta^dH_S(m-t+d)-\Delta^dH_S(m-t+1+d)=-\Delta^{d+1}H_S(m-t+1+d)<0$.
So $h_Y$ is of decreasing type.
\end{proof}

\begin{remark}
(i) From a result of Boij (see \cite{Boij}), it follows that for
$m\gg 0$ any aG divisor on an aCM scheme $S$ is linearly equivalent to $mH-K$, as in Theorem~\ref{lastprop}.\\
(ii) If we omit the hypothesis of irreducibility for $S$ in Theorem
\ref{lastprop},  we cannot say anything about the  irreducibility of
its twisted anti-canonical divisors.\\
(iii) Notice that $m \geq 2r-d$ implies that  the $h$-vectors of
the aG schemes $Y$ in Theorem
\ref{lastprop} have length $m+d$.\\
(iv) We were not able to compute an upper bound for the regularity of $\omega_S^{\vee}$
in terms
of invariants of $S$ such as its Betti numbers. However,
using a computer algebra system, such as CoCoA or
Macaulay2, it is possible to compute this regularity in concrete examples. In
Examples \ref{eg1} and \ref{eg2} we compute this bound for two concrete cases.\\
(v) Notice that it is in fact enough to take $m\ge\max\{2r-d,\alpha \}$, where $\alpha$ is the highest degree of a
minimal generator of $H^0_*(\omega_S^\vee)=Hom_S(K_S,S)$, $K_S$ the canonical module of $S$.
\end{remark}

Now we use the results of Section 2 and Theorem \ref{lastprop} to
obtain aG irreducible schemes as divisors on standard
determinantal schemes and to determine their Hilbert functions
 in terms  of the degree matrix.

\begin{corollary}\label{last.thm}Let $U=(u_{ij})$ be a degree matrix of
size $l\times (l+c-1)$. Let $n$ be an integer such that $1\leq c
\leq n-2$. Suppose that $u_{i,i+c}>0$. Then for $m\ge\max\{2r-n+c,\re (\omega_S^{\vee})\}$,
there exists an irreducible and reduced arithmetically Gorenstein
$Y\subset\PP ^n$ of codimension $c+1$ with the $h$-vector given by
\begin{enumerate}
\item for $t<r$,
$$h_Y(t)=\sum_{(V_1|V_2|\ldots)\subset U}(-1)^{1+\dim V_2+\ldots}
\bin{t+c-\tr(V_1|V_2|\ldots)}{c};$$
\item for $r\le t\le m-r+n-c$, $$h_Y(t)=
\frac{1}{c!}\sum_{(V_1|V_2|\ldots)\subset U}(-1)^{1+\dim V_2+\ldots+n}
(\tr(V_1|V_2|\ldots)-1)\cdot\ldots\cdot(\tr(V_1|V_2|\ldots)-c);$$
\item for $t\geq m-r+n-c$, $$h_Y(t)=
\sum_{(V_1|V_2|\ldots)\subset U}(-1)^{1+\dim
V_2+\ldots}\bin{m-t+n-\tr(V_1|V_2|\ldots)}{c}$$
\end{enumerate}
where $r={\tr(V_1^0 | V_2^0| \dots|V_c^0)}-c$, $V_1^0$ is the
submatrix of $U$ formed by the first $l$ columns and
$V_{k+1}^0=(u_{l,l+k})$ for $k=1, \dots c-1$.
\end{corollary}

\begin{proof}
We can choose a homogeneous matrix $A$ with degree matrix $U$ such
that it defines a standard determinantal scheme. From a theorem of
Trung (\cite{T}, Theorem 4.4 - notice that, since the polynomial
ring $R$ is normal, then normality carries on to the reduced,
irreducible lifting), there exists a normal, reduced, standard
determinantal scheme $S \subset \PP^n$ of dimension $d=n-c$ (so $d
\geq 2$) with associated degree matrix $U$. We remark that the aCM
scheme $S$ satisfies property $G_1$. Indeed, by Serre's Criterion
(see Theorem 11.5 of \cite{Eis}), $S$ satisfies property $G_1$
because it is irreducible and normal.

By Theorem \ref{lastprop}, $r=\re (S)-1$, and by Proposition \ref{hilbstd}, $\re (
S)={\tr(V_1^0|V_2^0| \dots |V_c^0)}-c+1$. Hence
$r={\tr(V_1^0|V_2^0| \dots |V_c^0)}-c$.

Next, we compute $\Delta^dH_S(t)$. By Proposition \ref{hilbstd},
we have that the $d$-th difference of $H_S$ is given by
\begin{align*}
\Delta ^dH_S(t)
%&=\sum _{-1\le i\le c-1}(-1)^{i+1}\Delta ^dH_{M_i}(t)\\
 &=\bin{t+n-d}{n-d}+\sum_{(V_1|V_2|\ldots)\subset U}(-1)^{1+\dim V_2+\ldots}
\bin{t+n-\tr(V_1|V_2|\ldots)-d}{n-d}\\
 &=\bin{t+c}{c}+\sum_{(V_1|V_2|\ldots)\subset U}(-1)^{1+\dim V_2+\ldots}
\bin{t+c-\tr(V_1|V_2|\ldots)}{c}.
\end{align*}

The Corollary follows  form Theorem \ref{lastprop}
and the expression of $\degr(S)$ obtained in Proposition
\ref{hilbstd}.
\end{proof}

Under stricter numerical conditions on the degree matrix $U$, we
also obtain Hilbert functions of smooth arithmetically Gorenstein
schemes.

\begin{corollary}\label{smooth} Let $U=(u_{ij})$ be a degree matrix
of a homogeneous $l\times (l+c-1)$-matrix of polynomials in
$k[x_0,\ldots ,x_n]$, $1 \leq c \leq n-2$. Assume that $n \leq
2c+1$, $u_{i,j} \neq 0$ for all $i,j$, $u_{1,k}>0$ if $k+
[\frac{n-c}{2}]+1-n \leq 0$, and $u_{k+[\frac{n-c}{2}]+1-n,k}>0$
if $k+ [\frac{n-c}{2}]+1-n > 0$. Then for $m\ge\max\{2r-n+c,\re (\omega_S^{\vee})\}$
there exists a
smooth, irreducible, reduced arithmetically Gorenstein subscheme
$Y\subset\PP ^n$ of codimension $c+1$ with the $h$-vector given by
\begin{enumerate}
\item for $t<r$, $$h_Y(t)=\sum_{(V_1|V_2|\ldots)\subset U}(-1)^{1+\dim V_2+\ldots}
\bin{t+c-\tr(V_1|V_2|\ldots)}{c};$$
\item for $r\le t\le m-r+n-c$, $$h_Y(t)=
\frac{1}{c!}\sum_{(V_1|V_2|\ldots)\subset U}(-1)^{1+\dim
V_2+\ldots+n}(\tr(V_1|V_2|\ldots)-1)\cdot\ldots\cdot(\tr(V_1|V_2|\ldots)-c);$$
\item for $t\geq m-r+n-c$, $$h_Y(t)=
\sum_{(V_1|V_2|\ldots)\subset U}(-1)^{1+\dim
V_2+\ldots}\bin{m-t+n-\tr(V_1|V_2|\ldots)}{c}$$
\end{enumerate}
where $r={\tr(V_1^0 | V_2^0| \dots|V_c^0)}-c$, $V_1^0$ is the
submatrix of $U$ formed by the first $l$ columns and
$V_{k+1}^0=(u_{l,l+k})$ for $k=1, \dots c-1$.
\end{corollary}

\begin{proof} A result of Ein \cite{Ein}, Theorem 2.6, ensures
that under the hypotheses of the Corollary, there exists a smooth standard
determinantal $S \subset \PP^n$ scheme with degree matrix $U$. Theorem
\ref{lastprop} applies to this $S$ giving a smooth aG scheme $Y$ (in the proof of
Theorem \ref{lastprop}, if $S$ is smooth then the Bertini Theorem gives
that $Y \in |mH-K|$ is also smooth). The computation of the Hilbert function of
$Y$ follows as in the proof of Corollary \ref{last.thm}.
\end{proof}

Here are two examples of $h$-vectors of irreducible aG schemes obtained using this technique.

\begin{example}\label{eg1}
Consider the irreducible and reduced standard determinantal scheme
$S \subset \PP^5$ associated to the matrix
$$\left( \begin{array}{cccc} x_0 & x_1 & x_2 & x_3 \\
x_2 & x_3 & x_4 & x_5 \end{array} \right)$$ It is a rational
normal scroll surface of $\PP^5$ of degree 4, whose $h$-vector is
$(1,3)$. Then the $h$-vector of a general $Y\in |mH-K|$ has form
$(1,4,4,\ldots, 4,1)$ and has length $m+2$. The Castelnuovo-Mumford regularity of the dual of the
canonical sheaf of the surface equals 5. Then $Y$ is irreducible for
$m \geq 5$.
\end{example}

Let $d$ be the degree of a reduced and irreducible subscheme of $\PP^n$ of codimension $3$. Then we may obtain for $Y$ as in
Corollary~\ref{last.thm} a Gorenstein, codimension $4$ $h$-vector containing an arbitrarily long constant
sequence of $d$'s in the middle.

\begin{example}\label{eg2}
Consider the degree matrix $$\left( \begin{array}{cccc} 2 & 2 & 2 & 1 \\
3& 3 & 3 & 2 \end{array} \right)$$
associated to the irreducible curve of $\PP^4$ whose defining matrix is
$$\left( \begin{array}{cccc} x_0^2  & x_1^2 & x_2^2 & x_3  \\
x_4^3  & x_0^3  & x_1^3  & x_3^2 \end{array} \right)$$ The curve
has a reduced irreducible lifting $S\subset\PP^5$, whose
$h$-vector is $(1,3,6,10,12,9,4,1)$. The $h$-vector of a general
divisor $Y$ on $S$ linearly equivalent to $mH-K$ has the form
$(1,4,10,20,32,41,45,46,\ldots,46,45,41,32,20,10,4,1)$, and has
length $m+2$. \newline
The Castelnuovo-Mumford regularity of the dual of the canonical sheaf
of the surface equals $8$. We also need $m\geq 2r-d=15$.
Then $Y$ is an irreducible codimension 4 arithmetically
Gorenstein scheme for $m\geq 15$.
\end{example}

\section*{Appendix}

We give here the proofs of the folklore results of Section
\ref{folklore}. For the notation, see Section \ref{folklore}.

\noindent {\it Proof of Proposition \ref{hilbstd}.}
 To start the computation of the
Hilbert function of a standard determinantal scheme $S\subset\PP
^n$, recall that a minimal free resolution of $R/I_S$ is given by
the Eagon-Northcott complex (see \cite{Eis}, Corollary A2.12):
\begin{align*}
M_*: 0 & \lra \bigwedge^{l+c-1} F\otimes
S_{c-1}(G)^*\otimes\bigwedge ^lG^*\lra\cdots\lra
\bigwedge ^{l+1}F\otimes S_{1}(G)^*\otimes\bigwedge ^lG^*\\
 & \lra \bigwedge ^{l}F\otimes S_{0}(G)^*\otimes\bigwedge ^lG^*\lra R\lra R/I_S\lra 0,
\end{align*}
where $\wedge (P), S(P),$ and $P^*$ mean the exterior algebra, the
symmetric algebra, and respectively, the dual of $P$ over $R$ for
any $R$-module $P$. Using the fact that $\wedge (P\oplus
P')=\wedge (P)\otimes\wedge (P')$ and $S(P\oplus P')=S(P)\otimes
S(P')$, we get for $i \geq 0$,
\begin{align*}
M_{i+1} &=\bigwedge^{l+i} F\otimes S_{i}(G)^*\otimes\bigwedge ^lG^*\\
 &=\bigwedge^{l+i}(\bigoplus _{1\le j\le l+c-1}R(a_j)) \otimes S_{i}
(\bigoplus _{1\le k\le l}R(b_k))^*\otimes\bigwedge ^l(\bigoplus _{1\le k\le l}R(b_k))^*\\
 &\cong (\bigoplus _{1\le j_1<\ldots <j_{l+i}\le l+c-1}R(a_{j_1}+\ldots+a_{j_{l+i}}))
\otimes(\bigoplus _{1\le k_1\le \ldots\le k_i\le l}R(-b_{k_1}-\ldots-b_{k_i}))\otimes\\
 & \otimes R(-b_1-\ldots-b_l)\\
 &\cong \bigoplus _{1\le j_1<\ldots <j_{l+i}\le l+c-1}^{1\le k_1\le \ldots\le
 k_i\le l}R(a_{j_1}+\ldots+a_{j_{l+i}}-b_{k_1}-\ldots-b_{k_i}-b_1-\ldots-b_l)\\
 &\cong \bigoplus
_{1 \le j_1 < \ldots < j_{l+i} \le l+c-1}^{1 \le k_1 \le \ldots
\le k_i \le
l}R(-u_{1,j_1}-\ldots-u_{l,j_l}-u_{k_1,j_{l+1}}-\ldots-u_{k_i,{j_{l+i}}}).
\end{align*}
By (\ref{pp}), after repeatedly replacing indices of the $u_{ij}$
by some permutations of them, we can write, for $i \ge 0$,
\begin{equation}\label{resol}
M_{i+1} \cong \bigoplus_{(V_1|V_2|\ldots)\subset U, \dim
V_2+\ldots=i}  R(-tr(V_1|V_2| \dots))
\end{equation}
Thus, for $i\ge 0$, we have
%\begin{align*}
%H_{M_{i+1}}(t) &=\sum_{1\le j_1<\ldots <j_{l+i}\le l+c-1}^{1\le
%k_1\le
%\ldots\le k_i\le l}\bin{t+n-u_{1,j_1}-\ldots-u_{l,j_l}-u_{k_1,j_{l+1}}-\ldots-u_{k_i,{j_{l+i}}}}{n}.\\
%\end{align*}
%I've put this sentence after the notation
\begin{align*}
H_{M_{i+1}}(t) &=\sum_{(V_1|V_2|\ldots)\subset U, \dim
V_2+\ldots=i}\bin{t+n-\tr(V_1|V_2|\ldots)}{n}.\\
\end{align*}
Therefore,
\begin{align*}
H_S(t) &=\sum _{0\le i\le c}(-1)^{i}H_{M_i}(t)\\
 &=\bin{t+n}{n}+\sum_{(V_1|V_2|\ldots)\subset U}(-1)^{1+\dim V_2+\ldots}\bin{t+n-\tr(V_1|V_2|\ldots)}{n}.
\end{align*}

At this point we may simplify. First, since $S$ has codimension
$c$, the coefficient of $t^i$ in $H_S(t)$ has to vanish for
$d+1\le i\le n$. Thus, we get
$$s_{n-i}(-1, \ldots,-n)+\sum_{(V_1|V_2|\ldots)\subset U}(-1)^{1+\dim V_2+\ldots}s_{n-i}
(\tr(V_1|V_2|\ldots)-1,\ldots,\tr(V_1|V_2|\ldots)-n)=0,$$ where
$s_j$ are the elementary symmetric functions in $n$ variables.
Similarly, since the $d$-th difference of the Hilbert function of
$S$
$$\Delta ^dH_S(t)
%&=\sum _{-1\le i\le c-1}(-1)^{i+1}\Delta ^dH_{M_i}(t)\\
% &=\bin{t+n-d}{n-d}+\sum_{(V_1|V_2|\ldots)\subset U}(-1)^{1+\dim V_2+\ldots}
%\bin{t+n-\tr(V_1|V_2|\ldots)-d}{n-d}\\
 =\bin{t+c}{c}+\sum_{(V_1|V_2|\ldots)\subset U}(-1)^{1+\dim V_2+\ldots}\bin{t+c-\tr(V_1|V_2|\ldots)}{c}
$$
is the Hilbert function of a zero-scheme in $\PP^c$, for $1\le
j\le c$ we have that
$$s_{c-j}'(-1, \dots,-c)+\sum_{(V_1|V_2|\ldots)\subset U}(-1)^{1+\dim V_2+\ldots}s_{c-j}'
(\tr(V_1|V_2|\ldots)-1,\ldots,\tr(V_1|V_2|\ldots)-c)=0,$$ where
$s_j'$ are the elementary symmetric functions in $c$ variables.

Let us call $r$ the length of the $h$-vector of $S$ or,
equivalently,
$$r:=\min\{i \mid \Delta^dH_S(i)=\deg S \}$$ (see the end of the
proof for an expression of $r$ in terms of the degree matrix
$U$). %Finally, if we define
%$$r'=\max\{\tr(V_1|V_2|\ldots)\ |\ (V_1|V_2|\ldots)\subset U\},$$
Then, by the vanishing formulas above, we have that for $t\ge r$,
$$\Delta ^{d}H_S(t)=\deg S$$
$$ =1+\frac{1}{c!}\sum_{(V_1|V_2|\ldots)\subset U}(-1)^{c+1+\dim V_2+\ldots}s_{c}'
(\tr(V_1|V_2|\ldots)-1,\ldots,\tr(V_1|V_2|\ldots)-c)$$
$$ =1+\frac{1}{c!}\sum_{(V_1|V_2|\ldots)\subset U}(-1)^{c+1+\dim V_2+\ldots}(\tr(V_1|V_2|\ldots)-1)
\cdot\ldots\cdot(\tr(V_1|V_2|\ldots)-c).$$

In order to obtain the expression about the regularity of $S$ we
look again at the Eagon-Northcott resolution of $R/I_S$. Notice
that, since $I_S$ is a perfect ideal of codimension $c$, the
$\max_{i,j} \{a_{i,j}-i \}$ will be achieved in the last free
module of a minimal free resolution of $I_S$ as an $R$-module.
Using the isomorphisms (\ref{resol}), Remark \ref{reg} and the
last observation, we get:

$$\re (S)=\max\{{\tr(V_1|V_2| \dots)}-(dim V_2+dim V_3+\ldots)\}=$$
$$=\max_i \{\max_
{\dim V_2+ \dots=i} \{{\tr(V_1|V_2| \dots)}-i \} \}=\max_{\dim
V_2+ \dots=c-1} \{{\tr(V_1|V_2| \dots)}-(c-1) \}.$$ From the way the $u_{i,j}$'s are ordered, we get that the maximum is
achieved by the combination of submatrices $W_{c-1}=(V_1^0|V_2^0|
\dots |V_{c}^0)$. The formula for the regularity follows.

Notice that the $r$ defined above satisfies $r=r(S)+d$ where
$r(S)$ is the index of regularity of $S$ (see Remark \ref{reg}).
Moreover, since $S$ is aCM, $r(S)=\re (S)-d-1$ (see Remark
\ref{reg}), so
$$r=\re (S)-1={\tr(V_1^0 | V_2^0| \dots|V_c^0)}-c.$$
\begin{flushright}$\Box$
\end{flushright}

\noindent {\it Proof of Proposition \ref{reduced}.}
 The idea of the proof is to see that there exists a reduced
standard determinantal scheme $S \subset \PP^m$, $m=\max \{ n,
2(l-1)+c-2 \}$, with degree matrix $U$. Then, taking $m-n$ general
hyperplane sections of it we get the desired reduced scheme $X
\subset  \PP^n$.

In what follows, whenever a claim involves a {\it general form} $G$
of degree $d$
in $R'=k[x_0,\ldots, x_m]$, it should be understood that the claim is true for
all $G$ outside a proper closed subset of the linear system $|\OO_{\PP^m}(d)|$.

We consider the matrix
$$A=\left(\begin{array}{cccccc}
G_1^1 & \cdots & G_1^c &  0  & 0& \cdots \\
0 & G_2^2 & \cdots & G_2^{c+1} & 0 &  \cdots \\
 & & \ddots  &  & \ddots & \\
0 & 0 & \cdots & G_l^l & \cdots & G_l^{l+c-1}
\end{array} \right)$$
where $G_i^j$ are general forms of degree $u_{i,j}$ in $k[x_0,
\dots, x_m]$, $m=\max \{ n, 2(l-1)+c-2 \}$. We are going to see
that the standard determinantal scheme $S \subset \PP^m$ defined
by the maximal minors of $A$ is a reduced scheme.

We proceed by induction on $c$.

If $c=1$, then $I_S=G_1^1 \cdot \ldots \cdot G_l^1$ and $S$ is the
union of $l$ general hypersurfaces, so $S$ is reduced.

If $c=2$, it follows from a result of Gaeta (\cite{G}). In this
case $S$ is a union of reduced complete intersections.

If $c \geq 3$, we will proceed by induction on $l$.

When $l=1$, $S$ is the complete intersection $(G_1^1, \dots,
G_1^c)$, so $S$ is reduced because $G_1^j$ are general.

When $l >1$, we claim that

\vspace*{2mm}

\underline{Claim 1:} $I_S=\bigcap\limits_{i=1}^l
((G_i^i)+I(B_i))+G_l^{l+c-1}I_Y$ where $I_Y$ is the ideal
generated by the maximal minors of the first $l-1$ rows and first
$l+c-2$ columns of $A$ and $I(B_i)$ is the ideal of maximal minors
of the following $(l-i+1) \times (l+c-2-i)$ submatrix of $A$:
$$B_i= \left(\begin{array}{ccccccc}
G_i^{i+1} & \cdots & G_i^{i+c-1} & 0 & \cdots & \cdots & 0 \\
G_{i+1}^{i+1} & \cdots & & G_{i+1}^{i+c} & 0 & \cdots & 0 \\
 & & & & & & \\
 & \ddots & & & & \ddots & \\
 & & & & & & \\
 & & & & & & \\
0 & \cdots & 0 & G^{l-1}_{l-1} & \cdots & \cdots & G_{l-1}^{l+c-2} \\
0 & \cdots & \cdots & 0 & G_l^{l} & \cdots & G_l^{l+c-2}
\end{array}\right).$$

\vspace*{2mm}

\emph{Proof of claim 1:} It is not difficult to check that
$I(B_{i+1})+(G_{i+1}^{i+1})\supseteq I(B_i)$. Hence, using the
modular law
  and that the
$G_{i}^i$ are general forms, we get that
$$\bigcap_{i=1}^l ((G_i^i)+I(B_i))+G_l^{l+c-1}I_Y=I(B_1)+ \sum_{i=1}^{l-1}
G_1^1 \cdot G_2^2 \cdot \ldots \cdot G_i^i \cdot I(B_{i+1})+G_1^1
\cdot G_2^2 \cdot \ldots \cdot G_l^l+G_l^{l+c-1}I_Y.$$

This last ideal is the ideal generated by the maximal minors of
$A$. Indeed, the maximal minors of $A$ that contain the last
column generate $G_l^{l+c-1}I_Y$. We restrict now to the minors
that do not contain the last column. Among them, the minors that
do not contain the first column generate $I(B_1)$; if we consider
the ones that contain the first column, then we must distinguish
between the ones that do not contain the second column, these
generate $G_1^1 I(B_2)$, and the ones that contain the second
column: for these we distinguish between the minors that do not
contain the third column, these generate $G_1^1  G_2^2 I(B_3)$,
and the ones that contain the third column, and so on. Hence
$$I(B_1)+ \sum_{i=1}^{l-1}
G_1^1 \cdot G_2^2 \cdot \ldots \cdot G_i^i \cdot I(B_{i+1})+G_1^1
\cdot G_2^2 \cdot \ldots \cdot G_l^l+G_l^{l+c-1}I_Y=I_S$$ and the
claim is proved.

\vspace*{3mm}

By induction on $l$, we know that $Y$ is a reduced standard
determinantal scheme of codimension $c$. Moreover, each $B_i$
defines a standard determinantal scheme $X_i$ of codimension
$c-2$.

\underline{Claim 2:} $B_i$ defines a reduced standard
determinantal scheme $X_i$.

\emph{Proof of claim 2:} %By induction hypothesis on $c$, $B_i$ defines a reduced standard
%determinantal scheme $X_i$ ($B_i$ is a deformation of the matrix
%that gives a reduced standard determinantal scheme in codimension
%$c-2$, which is obtained from $B_i$ when the subdiagonal is 0 and
%$G_{l-1}^{l+c-2}=0$).
By induction hypothesis we have that
$$C_i= \left(\begin{array}{ccccccc}
G_i^{i+1} & \cdots & G_i^{i+c-2} & 0 & \cdots & \cdots & 0 \\
0 & G_{i+1}^{i+1} & \cdots & G_{i+1}^{i+c-1} & 0 & \cdots & 0 \\
 & & \ddots & & & \ddots & \\
0 & \cdots & 0 & G^{l-1}_{l-1} & \cdots & G_{l-1}^{l+c-3} & 0 \\
0 & \cdots & \cdots & 0 & G_l^{l} & \cdots & G_l^{l+c-2}
\end{array}\right)$$
is associated to a reduced standard determinantal scheme of
codimension $c-2$. Moreover, if we denote by $I(B_i)$, $I(C_i)$
the ideals generated by the maximal minors of $B_i$, $C_i$
respectively, we have that
\begin{eqnarray*}I(C_i)+(G_{i+1}^{i+1},\ldots,G_l^l,G_i^{i+c-1},G_{i+1}^{i+c},
\ldots,G_{l-1}^{l+c-2}) =\\
I(B_i)+(G_{i+1}^{i+1},\ldots,G_l^l,G_i^{i+c-1},G_{i+1}^{i+c},
\ldots,G_{l-1}^{l+c-2}). \end{eqnarray*}

If we call $R'=k[x_0, \dots , x_m]$, we observe that
$G_{i+1}^{i+1},\ldots,G_l^l,G_i,^{i+c-1},\ldots,G_{l-1}^{l+c-2}$
is an $R'/I(C_i)$-regular sequence (since this is a regular
sequence in $R'$, $\dim{R'/I(C_i)}=m-c+3 \geq 2(l-1)+1$, and
$G_{i+1}^{i+1},\ldots, G_l^l,G_i,^{i+c-1},\ldots,G_{l-1}^{l+c-2}$
do not appear in $I(C_i)$). Hence, the ideal
$I(C_i)+(G_{i+1}^{i+1},\ldots,G_l^l,G_i,^{i+c-1},\ldots,G_{l-1}^{l+c-2})$
defines a reduced, standard determinantal scheme of codimension
$c-2$.

Then
$I(B_i)+(G_{i+1}^{i+1},\ldots,G_l^l,G_i,^{i+c-1},\ldots,G_{l-1}^{l+c-2})$
also defines a reduced, aCM scheme of codimension $c-2$. Since
$G_{i+1}^{i+1},\ldots,G_l^l, G_i,^{i+c-1},\ldots,G_{l-1}^{l+c-2}$
is a regular sequence modulo $I(B_i)$, we have that $X_i$ is a
reduced scheme of codimension $c-2$ (see \cite{CN}). This is the
end of the proof of Claim 2.

\vspace*{3mm}

Therefore, $$S=Y \cup \bigcup_{i=1}^l X_i\cap G_i^i \cap
G_l^{l+c-1},$$ and since $Y$ and $X_i\cap G_i^i \cap G_l^{l+c-1}$
are reduced, $S$ is also a reduced standard determinantal scheme
of codimension $c$ in $\PP^m$. Here, we call both $G_i^j$ the form and
the hypersurface defined by the form.

Now to obtain the desired reduced standard determinantal subscheme
$X \subset \PP^n$ we only need to take $m-n$ general hyperplane
sections of $S$. Notice that $m-n\geq m-c= \dim S$, so we are
taking general hyperplane sections of a reduced scheme and
reducibility is preserved.
\begin{flushright}$\Box$
\end{flushright}

\end{document}